\documentclass[letterpaper, 10 pt, conference]{ieeeconf}  

\usepackage{graphicx}          

\usepackage{epsfig} 
\usepackage{amsfonts}
\usepackage{amsmath,amssymb}
\usepackage{comment}
\usepackage{color}
\usepackage{cases}
\usepackage[subrefformat=parens]{subcaption}

\usepackage{standalone}
\usepackage{algpseudocode}
\usepackage{hyperref}
\usepackage{algorithm}
\usepackage{array}
\usepackage{textcomp}
\usepackage{stfloats}
\usepackage{url}
\usepackage{verbatim}
\usepackage{graphicx}

\usepackage{cite}

\usepackage{mathtools}				 

\usepackage{tikz}

\usepackage{bm}
\usepackage{here}
\usepackage{mathrsfs}

\usepackage{stmaryrd}

\makeatletter
\newcommand{\specialcell}[1]{\ifmeasuring@#1\else\omit$\displaystyle#1$\ignorespaces\fi}
\makeatother

\newtheorem{problem}{Problem}
\newtheorem{proposition}{Proposition}

\newtheorem{theorem}{Theorem}

\newtheorem{assumption}{Assumption}

\newcommand{\calE}{{\mathcal E}}    
    
\newcommand{\calG}{{\mathcal G}}

\newcommand{\calM}{{\mathcal M}}

\newcommand{\calR}{{\mathcal R}}

\newcommand{\scrB}{{\mathscr B}}

\newcommand{\scrG}{{\mathscr G}}

\newcommand{\bbR}{{\mathbb R}}

\newcommand{\bbZ}{{\mathbb Z}}

\newcommand{\rmd}{{\rm d}}

\newcommand{\rmh}{{\rm h}}

\newcommand{\rmH}{{\rm H}}

\newcommand{\sfd}{{\sf d}}
\newcommand{\sfe}{{\sf e}}

\newcommand{\bfone}{{\bf 1}}

\newcommand{\bmA}{{\boldsymbol A}}
\newcommand{\bmB}{{\boldsymbol B}}

\newcommand{\bmG}{{\boldsymbol G}}
\newcommand{\bmK}{{\boldsymbol K}}
\newcommand{\bmP}{{\boldsymbol P}}

\newcommand{\bmc}{{\boldsymbol c}}
\newcommand{\bmell}{{\boldsymbol \ell}}
\newcommand{\bmu}{{\boldsymbol u}}
\newcommand{\bmx}{{\boldsymbol x}}

\newcommand{\Ee}{{\rm e}}   

\newcommand{\bbra}[1]{\ensuremath{[\![#1]\!]} }  

\newcommand{\diagbox}{\boxbslash}

\newcommand{\mpc}{{\rm MPC}}
\newcommand{\tmp}{{\rm tmp}}

\newcommand{\calTone}{{\mathcal T}(\bfone_N/N)}

\newcommand{\wtilde}[1]{{\widetilde{#1}}}

\allowdisplaybreaks[1]

\IEEEoverridecommandlockouts                              
\title{\LARGE \bf
Entropic Model Predictive Optimal Transport for Underactuated Linear Systems
}

\author{Kaito Ito and Kenji Kashima
\thanks{
	\copyright 2023 IEEE. Personal use of this material is permitted. Permission from IEEE must be obtained for all other uses, in any current or future media, including reprinting/republishing this material for advertising or promotional purposes, creating new collective works, for resale or redistribution to servers or lists, or reuse of any copyrighted component of this work in other works.
}
\thanks{This work was supported in part by JSPS KAKENHI Grant Numbers JP21J14577, JP21H04875, and JST, ACT-X Grant Number JPMJAX2102.}
\thanks{K. Ito is with the School of Computing, Tokyo Institute of Technology,
Yokohama, Japan
        {\tt\small ka.ito@c.titech.ac.jp}}%
\thanks{K. Kashima is with the Graduate School of Informatics, Kyoto University,
Kyoto, Japan
        {\tt\small kk@i.kyoto-u.ac.jp}}%
}

\begin{document}

\maketitle
\thispagestyle{empty}
\pagestyle{empty}

\begin{abstract}
	This letter investigates dynamical optimal transport of underactuated linear systems over an infinite time horizon.
	In our previous work, we proposed to integrate model predictive control and the celebrated Sinkhorn algorithm to perform efficient dynamical transport of agents. However, the proposed method requires the invertibility of input matrices, which severely limits its applicability. To resolve this issue, we extend the method to (possibly underactuated) controllable linear systems. In addition, we ensure the convergence properties of the method for general controllable linear systems.
	The effectiveness of the proposed method is demonstrated by a numerical example.
\end{abstract}

\section{Introduction}\label{sec:intro}
The studies of large-scale systems composed of multiple agents are becoming increasingly important in view of their application such as sensor networks, smart grids, intelligent transportation systems, and systems biology.
Several topics have been investigated including formation control and synchronization. They can be expressed as the problem of stabilizing the distribution of agents to a desired distribution~\cite{Justh2003,Kuritz2019}. Especially when considering the efficiency of transporting agents, the above problem can be formulated as an optimal transport (OT) problem over dynamical systems~\cite{Bakshi2020,Krishnan2018}.

The dynamical OT problem requires determining where and how to transport each agent.
For example, when the distribution of agents and the target distribution are given by empirical distributions, where and how to transport the agents correspond to an optimal assignment problem and an optimal control~(OC) problem, respectively.
Efficient algorithms for solving (static) assignment problems have been developed such as the Hungarian algorithm~\cite{Kuhn1955} and the auction algorithm~\cite{Bertsekas1992}.
These methods have been successfully applied to multi-agent assignment problems; see e.g.,~\cite{Yu2014,Mosteo2017} and references therein.
However, when it comes to dynamical OT problems, the situation is more complicated. This is because the transport cost for the dynamical OT is obtained by solving an OC problem while in general, multi-agent assignment considers easily computable assignment (transport) costs such as distance-based costs. 
Especially when stabilizing agents efficiently, one needs to solve an infinite horizon OC problem, which is difficult to solve.

To circumvent this issue, our previous work~\cite{Ito2023sinkhorn} applied a model predictive control~(MPC) strategy to the dynamical OT. MPC achieves efficient control with a reasonable computational cost by solving a tractable finite horizon OC problem at each time instead of an infinite horizon OC problem~\cite{Mayne2014}.
In addition to OC problems, OT with MPC, which we call model predictive OT, recursively solves an assignment problem based on the current transport costs.
However, when the number of agents is large, solving an assignment problem at each sampling instant is computationally very expensive even with the Hungarian algorithm.
Then, \cite{Ito2023sinkhorn} resolved this issue by introducing entropy regularization. Entropy-regularized OT problems can be solved efficiently by the so-called Sinkhorn algorithm, which is highlighted by \cite{Cuturi2013}.
In view of this, we integrated the Sinkhorn algorithm and MPC to perform cost-effective dynamical transport.
The resulting method is called {\em Sinkhorn MPC} and reduces the computational burden for performing efficient transport. 
In addition, the convergence properties of Sinkhorn MPC have been revealed.
Despite its usefulness, Sinkhorn MPC requires the assumption that input matrices of agents following linear dynamics are invertible. This is a strong assumption and should be removed to extend the applicability of Sinkhorn MPC for example to underactuated systems.

In this letter, we generalize Sinkhorn MPC to input matrices which are possibly not full rank. This is done by making reasonable assumptions such as the controllability of the agents.
Moreover, we show the convergence properties and ultimate boundedness of Sinkhorn MPC for general linear systems.

\textit{Organization:}
The remainder of this letter is organized as follows. In Section~\ref{sec:previous}, we briefly introduce our previously proposed method and its properties. In Section~\ref{sec:global}, we extend Sinkhorn MPC to general input matrices and reveal its global convergence property. In Section~\ref{sec:local}, we show the ultimate boundedness and local asymptotic stability for the extended method. In Section~\ref{sec:example}, a numerical example illustrates the obtained results.
In Section~\ref{sec:conclusion}, we present our conclusions.

\textit{Notation:}
Let $ \bbR $ denote the set of real numbers.
The set of all positive (resp. nonnegative) vectors in $ \bbR^n $ is denoted by $ \bbR_{> 0}^n $ (resp. $ \bbR_{\ge 0}^n $). We use similar notations for the sets of all real matrices $ \bbR^{m\times n} $ and integers $ \bbZ $, respectively.
The set of integers $ \{1,\ldots,N\} $ is denoted by $ \bbra{N} $.
The Euclidean norm is denoted by $ \|\cdot \| $.
For a positive semidefinite matrix $ A $, denote $\|x\|_A := (x^\top A x)^{1/2}$.
The identity matrix of size $n$ is denoted by $I_n$ or $ I $ when its size is clear in the context. The matrix norm induced by the Euclidean norm is denoted by $ \| \cdot \|_2 $.
For vectors $ x_1,\ldots,x_N \in \bbR^n $, a collective vector $ [x_1^\top \ \cdots \ x_N^\top]^\top \in \bbR^{nN} $ is denoted by $ [x_1 ; \ \cdots \ ; x_N] $.
For $ \alpha = [\alpha_1 \ \cdots \ \alpha_N ]^\top \in \bbR^N $, the diagonal matrix with diagonal entries $ \{\alpha_i\}_{i=1}^N $ is denoted by $ \alpha^\diagbox $.
The element-wise division of $ a, b \in \bbR_{>0}^n $ is denoted by $ a \oslash b := [a_1 /b_1 \ \cdots \ a_n/b_n]^\top$.
The $ N $-dimensional vector of ones is denoted by $ \bfone_N $.
The gradient of a function $ f $ with respect to the variable $ x $ is denoted by $ \nabla_x f $.
For $ x, x' \in \bbR_{> 0}^n $, define an equivalence relation $ \sim $ on $ \bbR_{> 0}^n $ by $ x \sim x' $ if and only if $ \exists r> 0, x = r x' $.

\section{Brief Introduction of Sinkhorn MPC}\label{sec:previous}
In this section, we briefly introduce Sinkhorn MPC proposed in \cite{Ito2023sinkhorn} to efficiently solve a dynamical OT problem formulated as follows.
\begin{problem}
	\label{prob:main_continuous}
	Given initial and desired states $\{x_i^0\}_{i=1}^N, $ $ \{x_j^\sfd\}_{j=1}^N  \in (\mathbb{R}^{n})^N$,
	find control inputs $\{u_i\}_{i=1}^N$ and a permutation $ \sigma : \bbra{N} \rightarrow \bbra{N} $ that solve
	\begin{align}
		&\underset{\sigma }{\rm minimize} ~~ \sum_{i\in \bbra{N}} c_\infty^i (x_i^0, x_{\sigma(i)}^\sfd ) .  \label{eq:infinite_horizon_cost}
	\end{align}
	Here, the cost function $ c_\infty^i : \bbR^n \times \bbR^n \rightarrow \bbR $ is defined by
	\begin{align}
		& &&\hspace{-2cm} c_\infty^i (x_i^0, x_j^\sfd) := \min_{u_i} \ \int_0^\infty  \ell_i (x_i (t), u_i (t) ; x_j^\sfd) \rmd t \label{eq:value_infty}\\
		&\hspace{0.5cm} \text{subject to} && \dot{x}_i(t)  = A_i x_i (t)  + B_i u_i (t) , \label{eq:linear_dynamics_conti}\\
		& &&x_i(0) = x_{i}^0, \label{eq:constraint_initial} \\
		& &&\lim_{t\rightarrow \infty} x_i(t) =  x_j^\sfd,\label{eq:constraint_infty}
	\end{align}
	where $ x_i(t) \in \bbR^n $ denotes the state of the agent $ i $, and $ u_i(t)\in \bbR^m, A_i\in \bbR^{n\times n}, B_i \in \bbR^{n\times m} $.
	\hfill $ \diamondsuit $
\end{problem}

Here, a permutation $ \sigma $ determines the destination $ x_{\sigma(i)}^\sfd $ of each agent $ x_i $.
In order to satisfy the constraint~\eqref{eq:constraint_infty}, the system~\eqref{eq:linear_dynamics_conti} needs to admit a constant input $ u_i(t) \equiv \bar{u}_{ij} $ that makes $ x_i = x_j^\sfd $ an equilibrium of \eqref{eq:linear_dynamics_conti}. That is,
\begin{equation}\label{eq:constant_u}
        A_i x_j^\sfd + B_i \bar{u}_{ij} = 0.
\end{equation}
Note that the running cost $ \ell_i : \bbR^n \times \bbR^m \times \bbR^n \rightarrow \bbR $ is allowed to depend on the desired state $ x_j^\sfd $ so that at $ (x_j^\sfd,\bar{u}_{ij}) $, there is not a cost incurred, i.e., $ \ell_i(x_j^\sfd, \bar{u}_{ij}; x_j^\sfd) = 0 $.

Problem~\ref{prob:main_continuous} is challenging to solve due to the following two reasons:
First, in general, the infinite horizon OC problem~\eqref{eq:value_infty} is computationally intractable.
Second, even if $ c_\infty^i $ is available, it is still computationally expensive to solve the assignment problem~\eqref{eq:infinite_horizon_cost} when the number of agents $ N $ is large.

To avoid these issues, the previous work~\cite{Ito2023sinkhorn} proposed to utilize MPC and entropy regularization for OT. Instead of \eqref{eq:value_infty}, MPC recursively solves a tractable finite horizon OC problem in a receding horizon manner. The transport cost for the finite horizon problem with a prediction horizon $ T_\rmh > 0 $ is given by
\begin{align}
	c_{T_\rmh}^{i} (\check{x}_i, x_j^\sfd) := &\min_{u_i} \ \int_{0}^{T_\rmh} \ell_i (x_i(t), u_i (t); x_j^\sfd) \rmd t\label{eq:cost_finite_conti}\\
	&\text{subj. to  \eqref{eq:linear_dynamics_conti},}~~ x_i(0) = \check{x}_i, \  x_i (T_\rmh) = x_{j}^\sfd . \nonumber
\end{align}
The first input vector of the optimal input sequence of \eqref{eq:cost_finite_conti} is denoted by $ u_i^\mpc (\check{x}_i,x_j^\sfd) $.

Entropy regularization is useful to reduce the computational burden for solving static OT problems including the assignment problem~\eqref{eq:infinite_horizon_cost}~\cite{Cuturi2013,Kosowsky1994}. Let us consider the following problem, called an entropic OT problem:
\begin{equation}\label{prob:Sinkhorn}
        \underset{P\in \calTone}{\rm minimize} \ \sum_{i,j\in \bbra{N}} C_{ij} P_{ij} - \varepsilon \rmH(P),
\end{equation}
where 
\[
	\calTone := \bigl\{ P \in \bbR_{\ge0}^{N\times N} : P {\bfone}_N = P^\top \bfone_N = \bfone_N/N   \bigr\},
\]
$ C_{ij}:= c_\infty^i (x_i^{0},x_j^\sfd) $, $ \varepsilon > 0 $ is a regularization parameter, and the entropy of $ P $ is defined by $\rmH(P) := - \sum_{i,j} P_{ij} (\log (P_{ij}) - 1)$.
When $ \varepsilon = 0$, \eqref{prob:Sinkhorn} is a linear programming (LP), and it is known that there exists an optimal solution $ P^\sigma $ of the LP such that for the optimal permutation $ \sigma $ of \eqref{eq:infinite_horizon_cost}, it holds $ P_{ij}^\sigma = 1/N $ if $ j = \sigma(i) $, and $ 0 $, otherwise~\cite[Proposition~2.1]{Peyre2019}. 
Hence, $ \sigma $ can be recovered from $ P^\sigma $, which is called a permutation matrix for $ \sigma $, and in this sense, \eqref{prob:Sinkhorn} with $ \varepsilon= 0 $ is the {\em tight} LP relaxation of \eqref{eq:infinite_horizon_cost}.
In terms of mass transport, a matrix $ P \in \calTone $, which is referred to as a coupling matrix, represents a transport plan where $ P_{ij} $ describes the amount of mass flowing from $ x_i^{0} $ towards $ x_j^\sfd $.
In particular, $ P^\sigma $ expresses that all the mass $ 1/N $ at $ x_i^{0} $ is transported to $ x_{\sigma(i)}^\sfd $.

For large $ N $, solving the LP \eqref{prob:Sinkhorn} with $ \varepsilon = 0 $ is still computationally expensive.
An appealing feature of the entropy regularization ($ \varepsilon > 0 $) is that it admits an efficient algorithm.
Define the Gibbs kernel $ K $ associated with the cost $ C_{ij} $ as
\[
K = (K_{ij}) \in \bbR_{> 0}^{N\times N}, \ K_{ij} := \exp\left(  - C_{ij}/\varepsilon \right) .
\]
Then, a unique solution of \eqref{prob:Sinkhorn} has the form $ P^* =  (\alpha^*)^\diagbox K  (\beta^*)^\diagbox $ where the two scaling variables $ (\alpha^*,\beta^*) \in \bbR_{>0}^N \times \bbR_{>0}^N $ are determined by
\begin{equation}\label{eq:opt_scaling_condition}
\alpha^* = \bfone_N/N \oslash [K\beta^*], \ \beta^* = \bfone_N/N \oslash [K^\top \alpha^*] .
\end{equation}
The variables $ (\alpha^*, \beta^*) $ can be efficiently computed by the Sinkhorn algorithm:
\begin{align}
&\alpha[k+1] = \bfone_N/N \oslash [K\beta[k]], \ \beta[k] = \bfone_N/N \oslash \left[K^\top \alpha[k]\right], \nonumber\\
&\hspace{5.9cm} k\in \bbZ_{\ge 0}, \label{eq:Sinkhorn_algorithm}
\end{align}
where for any initial condition $ \alpha[0] = \alpha_0 \in \bbR_{>0}^N $, $\alpha[k+1]^\diagbox K \beta  [k]^\diagbox $ converges to $ P^* $ as $ k\rightarrow \infty $.
As $ \varepsilon $ goes to zero, the unique solution of \eqref{prob:Sinkhorn} converges to an optimal solution of \eqref{prob:Sinkhorn} without the regularization ($ \varepsilon = 0 $).
On the other hand, it is known that the convergence of the Sinkhorn algorithm deteriorates as $ \varepsilon \searrow 0 $.
Based on $ P^* $, the optimal destination $ x_{\sigma(i)}^\sfd $ of the agent $ x_i $ is typically approximated by the barycentric projection $ N \sum_{j=1}^N P_{ij}^* x_j^\sfd $~\cite[Remark~4.11]{Peyre2019}. Note that a permutation matrix $ P^\sigma $ satisfies $ x_{\sigma(i)}^\sfd =  N \sum_{j=1}^N P_{ij}^\sigma x_j^\sfd $.

To exploit the computational advantages of MPC and entropic OT, \cite{Ito2023sinkhorn} proposed to use the control law 
\begin{equation}\label{eq:entropic_mpc}
	u_i(t) = u_i^\mpc \bigl(x_i(t),x_i^\tmp (P^*(x(t))\bigr) ,
\end{equation}
where $ x(t) := [x_1(t);\cdots;x_N(t)] $, $ P^* (x) $ is the optimal solution of \eqref{prob:Sinkhorn} with $ C_{ij} = c_{T_\rmh}^i (x_i, x_j^\sfd), \ x = [x_1;\cdots;x_N] $,
and $ x_i^\tmp : \bbR_{\ge 0}^{N\times N} \rightarrow  \bbR^n $ determines a {\em temporary} target state of the $ i$th agent based on the coupling matrix $ P^*(x(t)) $ at each time $ t $.
This approach, which we call entropic model predictive OT (Ent-MPOT), substantially reduces the computational cost for performing efficient dynamical transport.

Moreover, considering the case in which only a small number of Sinkhorn iterations can be performed at each sampling instant, so that we cannot obtain a coupling close enough to $ P^*(x(t)) $, \cite{Ito2023sinkhorn} proposed to integrate MPC and the Sinkhorn algorithm. To explain this, we consider a time-discretized version of \eqref{eq:linear_dynamics_conti}, which is suitable for combining with the Sinkhorn iterations\footnote{Throughout this letter, we use bold symbols for discrete-time systems while we use italic letters for continuous-time systems.}:
\begin{equation}\label{eq:linear_dynamics_disc}
	\bmx_i [k+1]  = \bmA_i \bmx_i [k]  + \bmB_i \bmu_i [k], \ k\in \bbZ_{\ge 0} ,
\end{equation}
where $ (\bmA_i,\bmB_i) $ is obtained by e.g., a zero-order hold (ZOH) discretization of \eqref{eq:linear_dynamics_conti}.
Then, the cost function $ \bmc_{\tau_\rmh}^i $ with a finite horizon $ \tau_\rmh \in \bbZ_{>0} $ is defined by
\begin{align}
	\bmc_{\tau_\rmh}^i (\check{x}_i, x_j^\sfd) := &\min_{\bmu_i} \ \sum_{k=0}^{\tau_\rmh-1}  {\boldsymbol \ell}_i \left(\bmx_i [k], \bmu_i [k] ; x_j^\sfd \right) \label{eq:cost_disc}\\
	&\text{subj. to} ~~\eqref{eq:linear_dynamics_disc},\ \bmx_i[0] = \check{x}_{i},\ \bmx_i[\tau_\rmh] =  x_j^\sfd. \nonumber
\end{align}
Denote by $ \bmu_i^{\rm MPC} (\check{x}_i, x_j^\sfd) $ the first element of the optimal control sequence of the above problem. Let $ \bmx[k] := [\bmx_1[k];\cdots; \bmx_N[k]] $ and $ \bmP^* (x) $ be the optimal solution of \eqref{prob:Sinkhorn} with $ C_{ij} = \bmc_{\tau_\rmh}^i (x_i, x_j^\sfd) $.
Let $ S $ be the number of Sinkhorn iterations at each time $ k $.
Then Sinkhorn MPC integrating MPC and the Sinkhorn algorithm is given as follows.\\
{\bf Sinkhorn MPC:}
\begin{align}
&\bmx_i[k+1] = \bmA_i \bmx_i [k] + \bmB_i \bmu_i^{\rm MPC} \bigl(\bmx_i [k], x_i^\tmp \left(P[k] \right)  \bigr),\nonumber\\
&\hspace{6cm} \forall i \in \bbra{N}, \label{eq:SMPC_x} \\
&P[k] = \alpha\left[k,S+1\right]^\diagbox \bmK(\bmx[k]) \beta\left[k,S\right]^\diagbox, \label{eq:SMPC_P} \\
&\text{Sinkhorn iterations:} \nonumber\\
&\hspace{-0.2cm}\begin{cases}
\alpha\left[k,l+1\right] = \bfone_N/N \oslash \left[\bmK(\bmx[k])\beta[k,l] \right], \\
\beta[k,l] = \bfone_N/N \oslash \left[\bmK(\bmx[k])^\top \alpha[k,l] \right] , 
\end{cases} l\in \bbra{S} , \label{eq:SMPC_sink}\\
&\alpha[k+1,1] = \alpha[k,S+1],  \label{eq:SMPC_next} \\ 
&\bmx_i [0] = x_i^0, \ \alpha[0,1] = \alpha_0, \nonumber 
\end{align}
where
\begin{align*}
\bmK_{ij}(x) := \exp\biggl(- \frac{\bmc_{\tau_\rmh}^{i} (x_i, x_{j}^\sfd )}{\varepsilon} \biggr), \ x = [x_1;  \cdots  ;x_N] ,
\end{align*}
and the initial value $ \alpha_0 \in \bbR_{>0}^N $ is arbitrary.
\hfill $ \diamondsuit $

In addition to the computational efficiency of the proposed method, \cite{Ito2023sinkhorn} revealed its convergence properties and boundedness.
Specifically, when the running cost is quadratic, and $ x_i^\tmp $ is given by the barycentric projection, that is,
\begin{align}
        &\ell_i (x_i,u_i; x_j^\sfd) = \|u_i + B_i^{-1} A_i x_j^\sfd \|^2, \label{eq:quad_invertible_conti}\\
        &\bmell_i (x_i,u_i;x_j^\sfd) = \| u_i - \bmB_i^{-1} (x_j^\sfd - \bmA_i x_j^\sfd) \|^2, \\
        &x_i^\tmp (P) = \scrB_{i}(P,\{x_j^\sfd\}_{j=1}^N) :=  N \sum_{j=1}^{N} P_{ij} x_j^\sfd, \ P \in \bbR_{\ge 0}^{N\times N} , \label{eq:barycentric_target}
\end{align}
the following hold:
\begin{itemize}
        \item For any initial state, the solution of \eqref{eq:linear_dynamics_conti} under \eqref{eq:entropic_mpc} converges to the set of equilibrium points~\cite[Corollary~1]{Ito2023sinkhorn};
        \item The solution of \eqref{eq:SMPC_x} with \eqref{eq:SMPC_P}--\eqref{eq:SMPC_next} is ultimately bounded~\cite[Proposition~2]{Ito2023sinkhorn};
        \item For sufficiently small or large $ \varepsilon > 0 $, an equilibrium of \eqref{eq:SMPC_x}--\eqref{eq:SMPC_next} is locally asymptotically stable~\cite[Theorem~2]{Ito2023sinkhorn}.
\end{itemize}

However, the above results assume the invertibility of the input matrices $ \{B_i\} $ (or $ \{\bmB_i\} $) for two reasons. 
First, when $ B_i $ is not full row rank, that is, the $ i $th agent is underactuated, the existence of the inputs $ \{\bar{u}_{ij}\} $ in \eqref{eq:constant_u} does not ensure the existence of equilibria under Sinkhorn MPC. This is because $ \{x_j^\sfd\} $ is no longer an equilibrium under Sinkhorn MPC due to the regularization.
When $ B_i $ is full row rank, this issue does not arise because for any given state $ x_i^\sfe $, there always exists an input $ u_i $ that makes $ x_i^\sfe $ an equilibrium of \eqref{eq:linear_dynamics_conti}.
Second, we note that running costs $ \{\ell_i\} $ play a crucial role in ensuring the convergence of the proposed method like the conventional MPC~\cite{Mayne2000}. 
When $ B_i $ is not full column rank, there may be more than one constant input that makes a given state an equilibrium of \eqref{eq:linear_dynamics_conti}.
Then it is not trivial how the choice of such a constant input for designing the running cost $ \ell_i $ affects the dynamics under Sinkhorn MPC.

In the remainder of this letter, we remove the invertibility assumption. Specifically, by designing appropriate running costs which are quadratic in the control inputs, we reveal that even when the input matrices are not full rank, the barycentric projection~\eqref{eq:barycentric_target} resolves the above issues under reachability conditions.

\section{Global Convergence Property}\label{sec:global}
In this section, we deal with the continuous-time systems~\eqref{eq:linear_dynamics_conti} and consider the case where the number of Sinkhorn iterations at each time $ t $ is infinitely large.
In what follows, $ x_i^\tmp $ is given by the barycentric projection~\eqref{eq:barycentric_target}.
Instead of the invertibility of $ B_i $, we assume the following condition.
\begin{assumption}\label{ass:quadratic_cost}
	For all $ i\in \bbra{N} $, $ (A_i,B_i) $ is controllable.
	In addition, for all $ i,j \in \bbra{N} $, there exists $ \bar{u}_{ij} \in \bbR^m $ such that \eqref{eq:constant_u} holds.
\hfill$ \diamondsuit $
\end{assumption}

Hereafter, fix some $ \{\bar{u}_{ij}\} $ satisfying \eqref{eq:constant_u}.
Note that $ \bar{u}_i (P) := \scrB_{i} (P, \{\bar{u}_{ij}\}_{j=1}^N ) $ satisfies
\eqref{eq:constant_u} with $ \bar{u}_{ij} = \bar{u}_i (P), x_j^\sfd = x_i^\tmp (P) $.
This means that $ \bar{u}_i (P) $ is an equilibrium input that makes $ x_i = x_i^\tmp (P) $ an equilibrium of \eqref{eq:linear_dynamics_conti}.
Then we consider the following quadratic running cost depending on $ P $:
\begin{equation}\label{eq:quadratic}
	\ell_{i,P} \left(u_i \right) := \|u_i - \bar{u}_i (P)\|^2 , \ u_i \in \bbR^m, P\in \bbR_{\ge 0}^{N\times N} .
\end{equation}
Note that when $ B_i $ is invertible, the above cost coincides with \eqref{eq:quad_invertible_conti} with $ x_j^\sfd = x_i^\tmp (P) $.

Under the controllability of $ (A_i,B_i) $, the transport cost
\begin{align}
		&c_{T_\rmh,P}^{i} \left(\check{x}_i, x_i^\tmp (P)\right) := \min_{u_i} \ \int_{0}^{T_\rmh} \ell_{i,P} (u_i (t)) \rmd t \label{eq:prob}\\
		&\hspace{2cm} \text{subj. to  \eqref{eq:linear_dynamics_conti},}~~ x_i(0) = \check{x}_i, \  x_i (T_\rmh) = x_i^\tmp (P) , \nonumber
\end{align}
and the control law for MPC associated with \eqref{eq:prob} can be written as follows~\cite[Section~3.3,~pp.~138--140]{Lewis2012}:
\begin{align}
	&c_{T_\rmh,P}^{i}\left(\check{x}_i, x_i^\tmp (P)\right)  =\| \check{x}_i - x_i^\tmp (P) \|_{\scrG_i}^2, \label{eq:quad_value_func_conti}\\
	&u_{i,P}^{\rm MPC} \left(\check{x}_i, x_i^\tmp (P)\right) = - B_i^\top \scrG_i (\check{x}_i - x_i^\tmp (P) ) + \bar{u}_i (P) , \nonumber \\
	& \specialcell{\hfill \forall i \in \bbra{N}, \ \check{x}_i \in \bbR^{n}, \ P\in \bbR_{\ge 0}^{N\times N}}, \nonumber \\
	&\scrG_i := \biggl(\int_0^{T_\rmh} \Ee^{-A_i t} B_i B_i^\top  \Ee^{-A_i^\top t} \rmd t\biggr)^{-1} .\nonumber
\end{align}
Note that for a permutation matrix $ P^\sigma $, we have $ c_{T_\rmh,P^\sigma}^{i}(x_i, x_i^\tmp (P^\sigma))  = c_{T_\rmh,P^\sigma}^{i}(x_i, x_{\sigma(i)}^\sfd)=\| x_i - x_{\sigma(i)}^\sfd \|_{\scrG_i}^2 $.

Now, as a generalization of Ent-MPOT, we propose to use the controller
\begin{equation}
	u_i(t) = u_{i,P^*(x(t))}^{\rm MPC} \left(x_i(t), x_i^\tmp (P^*(x(t)))\right),
\end{equation}
where $ P^* (x)$ is the optimal solution of \eqref{prob:Sinkhorn} with $ C_{ij}= \|x_i - x_j^\sfd \|_{\scrG_i}^2 $.
	Then, the dynamics~\eqref{eq:linear_dynamics_conti} is written as
	\begin{align}
		\hspace{-0.22cm}\dot{x}_i (t) = (A_i - B_iB_i^\top \scrG_i) \biggl(x_i(t) - N \sum_{j=1}^N P_{ij}^*(x(t))x_j^\sfd \biggr) , \label{eq:lq_mpc_dynamics}
	\end{align}
	where we used the relationship $ B_i \bar{u}_i (P) = -A_i x_i^\tmp (P) $.
	We emphasize that \eqref{eq:lq_mpc_dynamics} no longer depends on the choice of $ \{\bar{u}_{ij}\} $.
Let us state the convergence result for \eqref{eq:lq_mpc_dynamics}.
The proof is given in Appendix~\ref{app:lasalle}.
\begin{theorem}\label{thm:convergence}
        Suppose that Assumption~\ref{ass:quadratic_cost} holds.
        Let
        \begin{align}
        \calM' &:= \biggl\{x\in \bbR^{nN} : B_i^\top \scrG_i \biggl(x_i - N\sum_{j=1}^N P_{ij}^*(x)x_j^\sfd \biggr) = 0, \nonumber\\
        & B_i^\top \Ee^{-A_i^\top T_\rmh} \scrG_i \biggl(x_i - N\sum_{j=1}^N P_{ij}^*(x)x_j^\sfd \biggr) = 0, \ \forall i\in \bbra{N} \biggr\} . \nonumber\label{eq:stationary_set}
        \end{align}
        Then, for any initial state $ x(0)  \in \bbR^{nN} $, $ x(t) $ following \eqref{eq:lq_mpc_dynamics} converges to the largest invariant set in $ \calM' $ for the dynamics
        \begin{equation}\label{eq:invariant_system}
                \dot{x}_i (t)= A_i \biggl(x_i(t) - N\sum_{j=1}^N P_{ij}^* (x(t)) x_j^\sfd \biggr),  \ i\in \bbra{N} . 
        \end{equation}
        \hfill $ \diamondsuit $
\end{theorem}

It is known that $ A_i - B_iB_i^\top \scrG_i $ is a Hurwitz matrix~\cite{Kleinman1970}, and thus the set of all equilibria of \eqref{eq:lq_mpc_dynamics} is given by $ \calR := \{x^\sfe \in \bbR^{nN} : x_i^\sfe  = N \sum_j P_{ij}^* (x^\sfe)x_j^\sfd , \forall i \in \bbra{N} \} $.
 The largest invariant set in $ \calM' $ contains $ \calR $. Especially when $ B_i $ is invertible for all $ i\in \bbra{N} $, it is obvious that $ \calR $ itself is the largest invariant set.
As mentioned in Section~\ref{sec:previous}, for any permutation $ \sigma $, $ x_\sigma^\sfd := [x_{\sigma(1)}^\sfd;\cdots;x_{\sigma(N)}^\sfd] \not\in \calR $ because $ P_{ij}^* (x) > 0 $ for any $ x\in \bbR^{nN} $. Similar to \cite[Lemma~1]{Ito2023sinkhorn}, it can be shown that for any $ \sigma $, there exists an equilibrium $ x^\sfe \in \calR $ which converges to the original target $ x_\sigma^\sfd $ as $ \varepsilon \searrow 0 $.

\section{Ultimate Boundedness and Local Asymptotic Stability}\label{sec:local}
Next, we explain that even for a finite number of Sinkhorn iterations, the boundedness for a generalized Sinkhorn MPC always holds, and in addition, a convergence result holds in a local sense.
We assume the following condition corresponding to Assumption~\ref{ass:quadratic_cost}.
\begin{assumption}\label{ass:quadratic_cost_disc}
        For all $ i\in \bbra{N} $, the reachability Gramian
        \begin{equation}
                \bmG_{i,\tau_\rmh} := \sum_{k=0}^{\tau_\rmh -1} \bmA_i^k \bmB_i \bmB_i^\top (\bmA_i^\top)^k
        \end{equation}
        is invertible.
        In addition, for all $ i,j \in \bbra{N} $, there exists $ \bar{\bmu}_{ij} \in \bbR^m $ such that
\begin{equation}\label{eq:constant_u_disc}
        \bmA_i x_j^\sfd + \bmB_i \bar{\bmu}_{ij} = x_j^\sfd 
\end{equation}
holds.
\hfill$ \diamondsuit $
\end{assumption}
Note that if \eqref{eq:linear_dynamics_disc}	is obtained by ZOH of a controllable system, there exists $ \tau_\rmh \in \bbZ_{> 0} $ such that $ \bmG_{i,\tau_\rmh} $ is invertible except for pathological cases~\cite[Theorem 3.2.1]{Chen1995}.
We emphasize that even if \eqref{eq:linear_dynamics_disc} is not a discretization of a continuous-time system and originally evolves in discrete time, the extended method and the results in this section can be applied to \eqref{eq:linear_dynamics_disc}.
This means that our idea can also be employed for discrete-time dynamical OT problems.

Similar to \eqref{eq:quadratic}, we consider the following quadratic cost:
\begin{align}
        \bmell_{i,P} (\bmu_i) &:= \| \bmu_i - \bar{\bmu}_i (P) \|^2,  \bmu_i \in \bbR^m, P \in \bbR_{\ge 0}^{N\times N} , \label{eq:general_quad_cost_disc}\\
        \bar{\bmu}_i(P) &:= \scrB_{i}(P,\{\bar{\bmu}_{ij}\}_{j=1}^N) .
\end{align}
The transport cost and the control law associated with \eqref{eq:general_quad_cost_disc} are given as follows~\cite[Section~2.2, pp.~37--39]{Lewis2012}:
\begin{align}
        \bmc_{\tau_\rmh,P}^{i} \left(x_i, x_i^\tmp(P)\right) &:=\| x_i - x_i^\tmp(P) \|_{\calG_i}^2, \label{eq:general_quad_value_func_disc}\\
        \bmu_{i,P}^{\rm MPC} \left(x_i,x_i^\tmp(P)\right) &:= - \bmB_i^\top (\bmA_i^\top)^{\tau_\rmh-1} \bmG_{i,\tau_\rmh}^{-1} \bmA_i^{\tau_\rmh} \nonumber\\
        &\quad\times (x_i - x_i^\tmp(P) ) + \bar{\bmu}_i (P) , \label{eq:general_mpc_quad_disc}
\end{align}
where $\calG_{i} := (\bmA_i^{\tau_\rmh})^\top \bmG_{i,\tau_\rmh}^{-1}  \bmA_i^{\tau_\rmh}$, $ i \in \bbra{N}, \ x_i \in \bbR^{n}, \ P\in \bbR_{\ge 0}^{N\times N} $.
Then, we propose a generalized Sinkhorn MPC as 
\begin{equation}\label{eq:general_sinkhorn_mpc}
	\bmu_i [k] = \bmu_{i,P[k]}^{\rm MPC} \left(\bmx_i[k],x_i^\tmp(P[k])\right) 
\end{equation}
whose $ P[k] $ is obtained by \eqref{eq:SMPC_P} with $ \bmK_{ij}(x) = \exp(-\|x_i-x_j^\sfd\|_{\calG_i}^2) / \varepsilon $.

For simplicity we consider only the case where just one Sinkhorn iteration is performed at each time, i.e., $ S = 1 $. Similar arguments in this section apply to the case where more iterations are performed. Then the dynamics under the generalized Sinkhorn MPC~\eqref{eq:general_sinkhorn_mpc} can be written as follows:
\begin{align}
  \hspace{-0.5cm}\bmx_i [k+1] &= \bar{\bmA}_i \bmx_i[k] + (I - \bar{\bmA}_i) x_i^\tmp \bigl(\wtilde{P}(\bmx[k],\beta[k])\bigr), \label{eq:xi_system}\\
	\beta[k+1] &= f (\bmx[k+1], \beta[k]) , \label{eq:xbeta_global}
\end{align}
where $ \bar{\bmA}_i := \bmA_i - \bmB_i \bmB_i^\top (\bmA_i^\top)^{\tau_\rmh -1} \bmG_{i,\tau_\rmh}^{-1} \bmA_i^{\tau_\rmh} $ and
\begin{align*}
				f (x,\beta) &:= \bfone_N/N \oslash \left[ \bmK(x)^\top (\bfone_N/N \oslash [\bmK(x) \beta]) \right], \\
        \wtilde{P}(x,\beta) &:=  \left( \bfone_N/N \oslash [\bmK(x)\beta]  \right)^\diagbox \bmK(x)  \beta^\diagbox  
\end{align*}
for $ (x,\beta) \in \bbR^{nN} \times \bbR_{>0}^N $.
In what follows, we regard $ \beta[\cdot] $ as a trajectory in the projective cone $ \bbR_{>0}^N / {\sim} $ because $ \alpha^*,\beta^* $ are only defined up to a multiplicative constant (see the Notation in Section~\ref{sec:intro} for $ \sim $).
Since $ \bmG_{i,\tau_\rmh} $ is invertible by Assumption~\ref{ass:quadratic_cost_disc}, $ \bar{\bmA}_i $ is stable, i.e., the spectral radius $ \rho_i $ of $ \bar{\bmA}_i $ satisfies $ \rho_i < 1 $~\cite[Theorem~1]{Selbuz1988}.
In addition, since $ \wtilde{P}(\bmx[k],\beta[k]) \bfone_N = \bfone_N /N $, the barycentric projection~\eqref{eq:barycentric_target} satisfies the following boundedness:
\begin{equation}\label{eq:barycentric_bound}
	\| x_i^\tmp(\wtilde{P}(\bmx[k],\beta[k])) \| \le \bar{r} := \max_{j\in \bbra{N}} \|x_j^\sfd\| .
\end{equation}
Therefore, \eqref{eq:xi_system} can be seen as a stable system whose input $ (I - \bar{\bmA}_i) x_i^\tmp (\wtilde{P}(\bmx[k],\beta[k])) $ is bounded. Then, by the same proof as in that of \cite[Proposition~2]{Ito2023sinkhorn}, we obtain the following.

\begin{proposition}\label{prop:bounded}
	Suppose that Assumption~\ref{ass:quadratic_cost_disc} holds.
	Then, for any $ \delta > 0, \{x_i^0\}_i $, and $ \{\nu_i\}_i $ satisfying $ \nu_i > 0, \rho_i + \nu_i < 1, \forall i\in \bbra{N} $, there exist $ \kappa_i(\nu_i) > 0, i\in \bbra{N} $ and $ \tau(\delta, \{x_i^0\},\{\nu_i \} ) \in \bbZ_{>0} $ such that the solution $ \bmx[k] = [\bmx_1[k];\cdots;\bmx_N[k]] $ of \eqref{eq:xi_system},~\eqref{eq:xbeta_global} satisfies
	\begin{equation}\label{ineq:ultimate_bound}
		\|\bmx_i [k] \| < \delta  +  \frac{ \kappa_i \bar{r} \|I- \bar{\bmA}_i\|_2}{1 - (\rho_i + \nu_i ) }  , \ \forall k \ge \tau, \ \forall i \in \bbra{N}. \tag*{$ \diamondsuit $}
	\end{equation}
\end{proposition}

Next in order to state the convergence result for \eqref{eq:xi_system},~\eqref{eq:xbeta_global}, we introduce some notation.
Note that equilibria $ (x^\sfe(\varepsilon), \beta^\sfe (\varepsilon)) \in \bbR^{nN} \times (\bbR_{>0}^N/ {\sim}) $ of \eqref{eq:xi_system},~\eqref{eq:xbeta_global} depend on $ \varepsilon $. 
By \cite[Lemma~1]{Ito2023sinkhorn}, for any permutation $ \sigma $, there exists an equilibrium $ (x^\sfe(\varepsilon), \beta^\sfe (\varepsilon)) $ of \eqref{eq:xi_system},~\eqref{eq:xbeta_global} such that $ x^\sfe (\varepsilon) $ and $ P^*(x^\sfe(\varepsilon)) $ converge exponentially to $ x_\sigma^\sfd $ and $ P^\sigma $, respectively, as $ \varepsilon \searrow 0 $. Denote by $ {\rm Exp}(\sigma) $ the set of all equilibria $ (x^\sfe(\cdot), \beta^\sfe (\cdot)) $ of \eqref{eq:xi_system},~\eqref{eq:xbeta_global} having the above exponential convergence property for a permutation $ \sigma $.
Now, we are ready to state the local convergence result.
Noting that the dynamics~\eqref{eq:xi_system},~\eqref{eq:xbeta_global} are exactly the same as (56), (57) in \cite{Ito2023sinkhorn}, we can apply the same proof as in \cite[Theorem~2]{Ito2023sinkhorn}.
\begin{theorem}\label{thm:local_convergence}
        Suppose that Assumption~\ref{ass:quadratic_cost_disc} holds, and for all $ i\in \bbra{N} $, $ \bmA_i $ is invertible. Also, assume that for some $ \varepsilon' > 0 $, $ (x^\sfe (\varepsilon'), \beta^\sfe(\varepsilon')) \in \bbR^{nN} \times (\bbR_{>0}^N / {\sim}) $ is an isolated equilibrium of \eqref{eq:xi_system},~\eqref{eq:xbeta_global}.
        Then the following hold:
        \begin{itemize}
                \item[(i)] For sufficiently large $ \varepsilon > 0 $, $ (x^\sfe(\varepsilon), \beta^\sfe(\varepsilon)) $ is locally asymptotically stable.
		\item[(ii)]  
		Assume further that $ x_i^\sfd \neq x_j^\sfd $ for all $(i,j), \ i\neq j $, and $ (x^\sfe (\cdot), \beta^\sfe(\cdot)) \in {\rm Exp}(\sigma) $ for some permutation $ \sigma $.
		Then, for sufficiently small $ \varepsilon >0$, $ (x^\sfe (\varepsilon), \beta^\sfe(\varepsilon)) $ is locally asymptotically stable.
                \hfill $ \diamondsuit $
        \end{itemize}
\end{theorem}

\section{Numerical Example}\label{sec:example}
\begin{figure}[t]
	\begin{minipage}[b]{1.0\linewidth}
		\centering
		\includegraphics[keepaspectratio, scale=0.29]
		{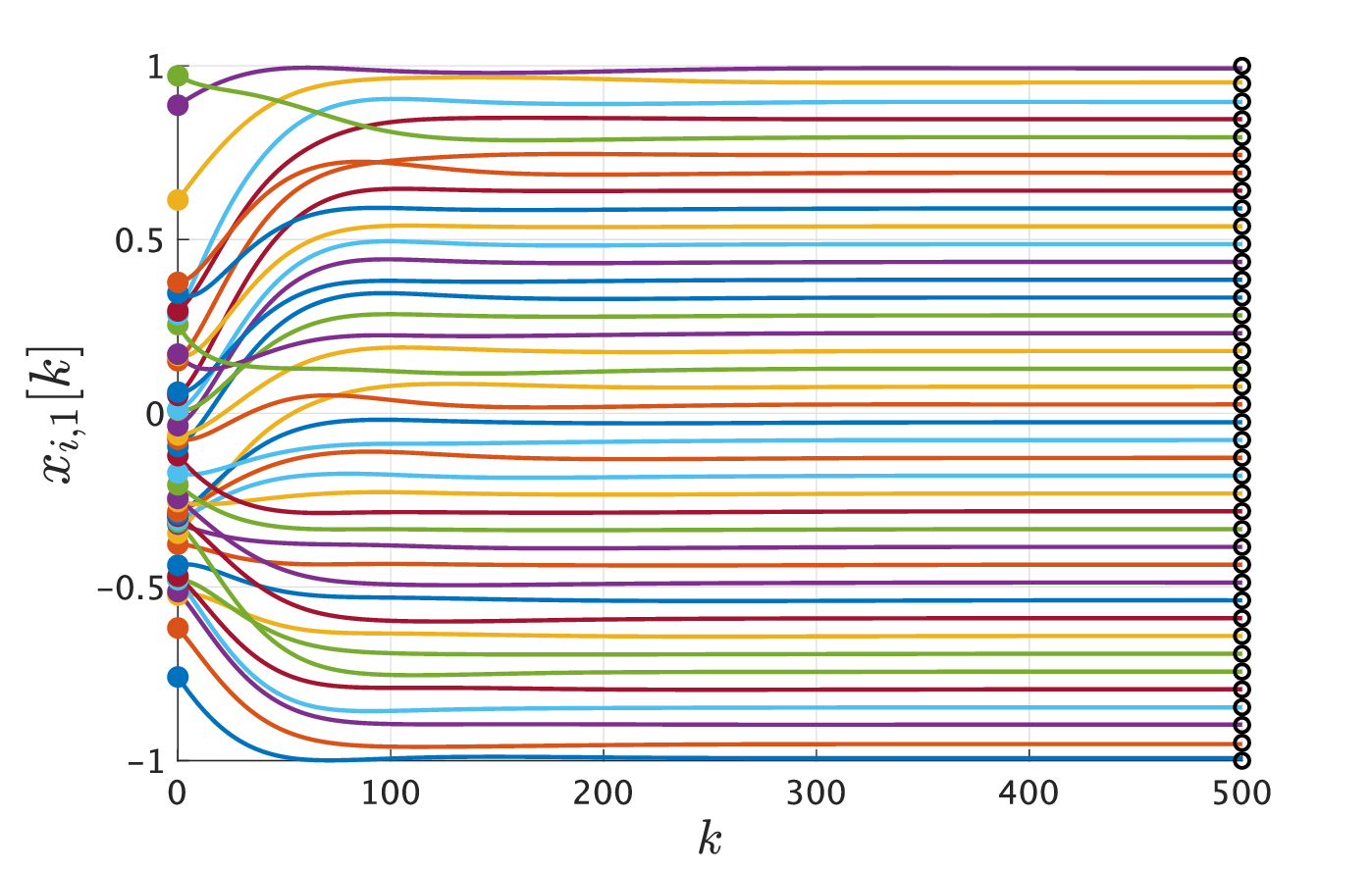}
	\end{minipage}
        \\
	\begin{minipage}[b]{1.0\linewidth}
		\centering
		\includegraphics[keepaspectratio, scale=0.29]
		{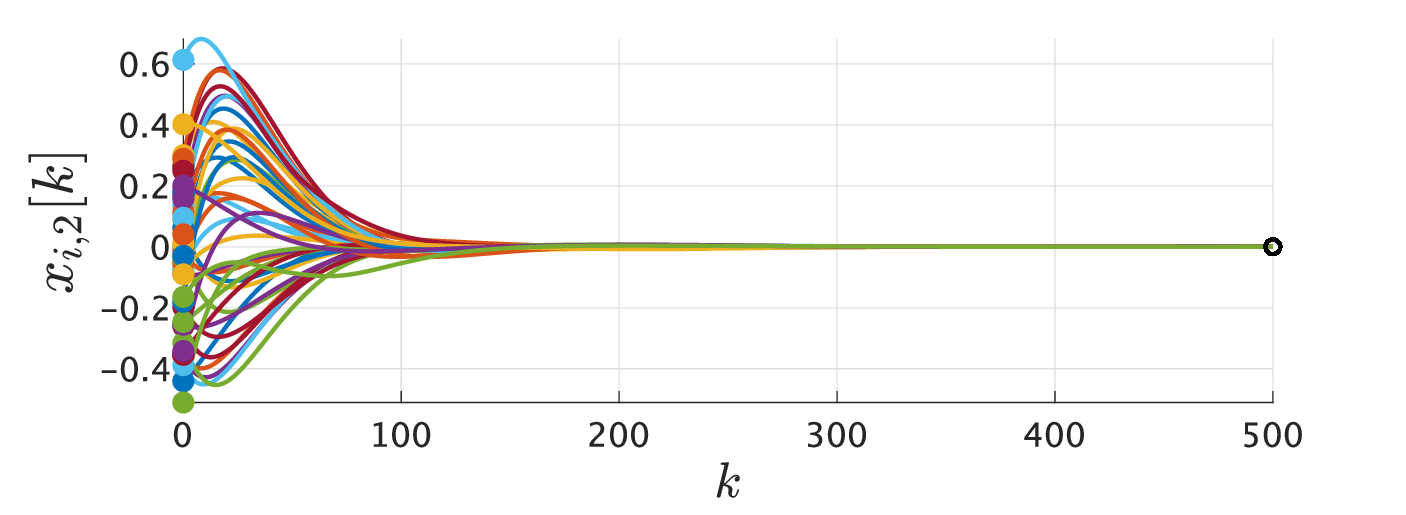}
	\end{minipage}
	\caption{Trajectories $ \bmx_i[k] = [\bmx_{i,1}[k] \ \bmx_{i,2}[k]]^\top $ of 40 agents for \eqref{eq:integrator_disc} with $ S = \infty $ (solid), initial states (filled circles), and desired states (black circles).}\label{fig:trajectory}
\end{figure}

\begin{figure}[t]
	\begin{minipage}[b]{1.0\linewidth}
		\centering
		\includegraphics[keepaspectratio, scale=0.30]
		{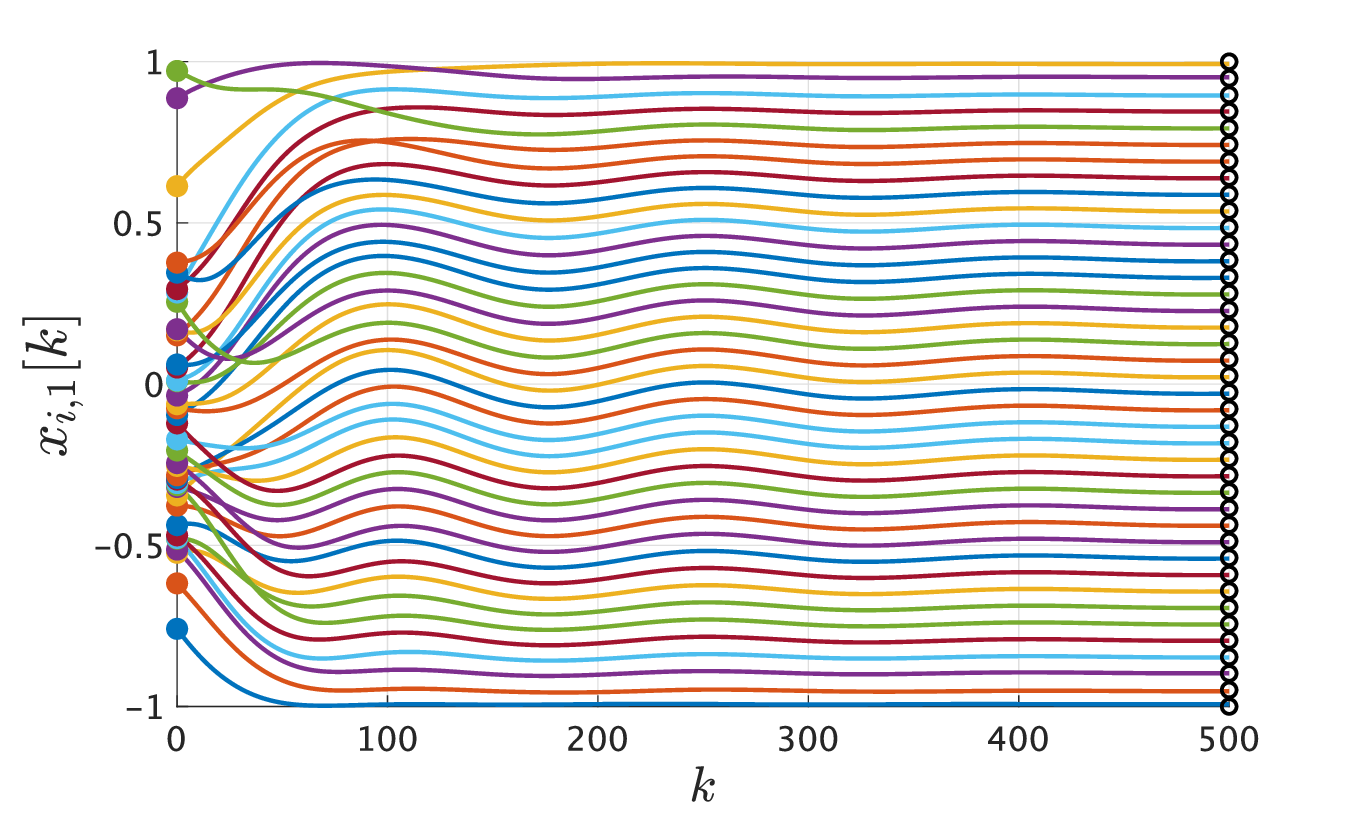}
	\end{minipage}
        \\
	\begin{minipage}[b]{1.0\linewidth}
		\centering
		\includegraphics[keepaspectratio, scale=0.30]
		{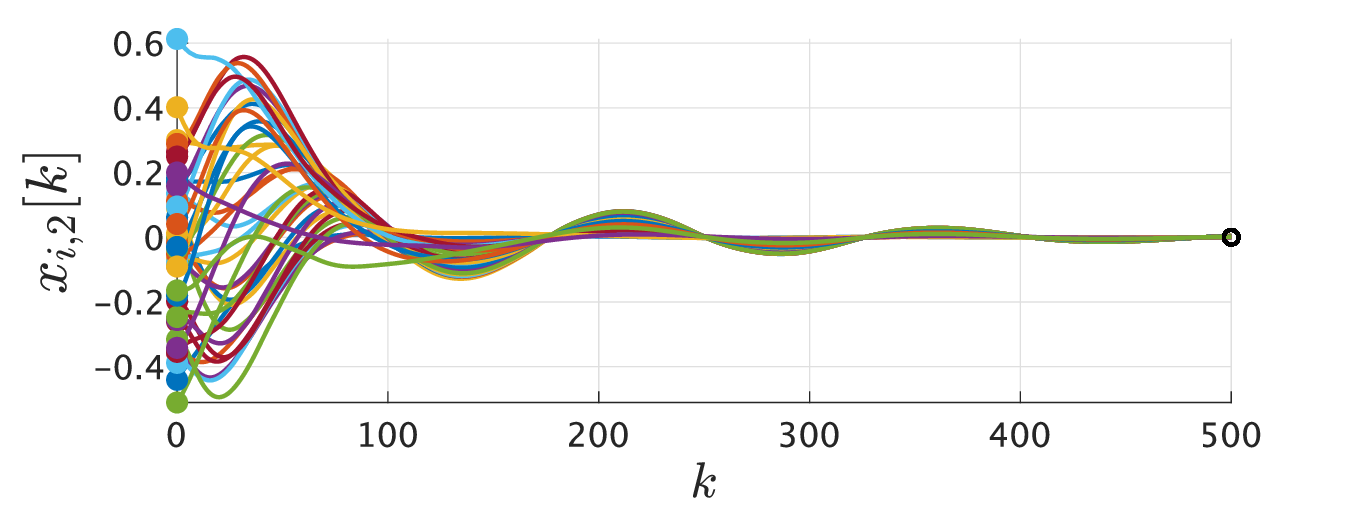}
	\end{minipage}
	\caption{Trajectories $ \bmx_i[k] = [\bmx_{i,1}[k] \ \bmx_{i,2}[k]]^\top $ of 40 agents for \eqref{eq:integrator_disc} with $ S = 20 $ (solid), initial states (filled circles), and desired states (black circles).}\label{fig:trajectory_finite}
\end{figure}

In this section, we illustrate how the generalized Sinkhorn MPC works via a numerical example.
We consider double integrator systems given by
\begin{equation}\label{eq:integrator_conti}
	A_i = 
	\begin{bmatrix}
		0 & 1\\
		0 & 0
	\end{bmatrix}, \
	B_i = 
        \begin{bmatrix}
                0\\1
        \end{bmatrix}, 
         \ \forall i\in \bbra{N} .
\end{equation}
Applying ZOH with the sampling period $ 0.02 $ to \eqref{eq:integrator_conti} yields
\begin{equation}\label{eq:integrator_disc}
	\bmA_i = 
	\begin{bmatrix}
		1 & 0.02\\
		0 & 1
	\end{bmatrix}, \
	\bmB_i = 
        \begin{bmatrix}
                0.0002\\0.02
        \end{bmatrix}, 
         \ \forall i\in \bbra{N} .
\end{equation}
The desired states $ \{x_j^\sfd \} $ are set to $ x_j^\sfd = [x_{j,1}^\sfd \ 0]^\top $ for some $ x_{j,1}^\sfd \in \bbR $. Then, $ \bar{u}_{ij} = 0 $ and $\bar{\bmu}_{ij} = 0 $ satisfy \eqref{eq:constant_u} and \eqref{eq:constant_u_disc}, respectively.
The regularization parameter and the prediction horizon are set to $ \varepsilon = 0.7, \tau_\rmh = 50 $.
Note that the trade-off in determining $ \varepsilon $ between the deviation of equilibria from the original target distribution and the transient behavior of the agents has already been discussed for the previous method in~\cite{Ito2023sinkhorn}, and the same argument holds for the extended method.
Fig.~\ref{fig:trajectory} depicts the trajectories $ \{\bmx_i[k] \}_i $ of \eqref{eq:linear_dynamics_disc} driven by \eqref{eq:general_sinkhorn_mpc} where Sinkhorn iterations are performed until the convergence criterion~\cite[Remark~4.14]{Peyre2019} is achieved. That is, $ \{x_i(t)\}_i $ following \eqref{eq:lq_mpc_dynamics} are well approximated by $ \{\bmx_i[k] \}_i $.
As expected from Theorem~\ref{thm:convergence}, the agents are successfully transferred to states that are sufficiently close to the desired states.
Next, Fig.~\ref{fig:trajectory_finite} shows $ \{\bmx_i[k]\}_i $ under Sinkhorn MPC with $ S = 20 $. Since the number of iterations $ S $ at each time is too small for $ P[k] $ not to be close enough to $ P^*(\bmx[k]) $ when $ k $ is small, the oscillations of the agents are observed.
Even in this case, the ultimate boundedness of $ \{\bmx_i[k]\}_i $ is ensured by Proposition~\ref{prop:bounded}.
Moreover, as can be seen, they still converge close to the desired states.

To see how the Sinkhorn algorithm combined with MPC affects control performance, we compare the accumulated cost $ \sum_{i,k} \Delta t \|\bmu_i[k]\|^2 ,  \Delta t = 0.02$ for the proposed method and MPC without the entropy regularization where $ P[k] $ is obtained by solving LP~\eqref{prob:Sinkhorn} with $ \varepsilon = 0 , C_{ij} = \|\bmx_i[k] - x_j^\sfd\|_{\calG_i} $.
Then the cost for the unregularized MPC is $ 11.56 $ while the cost for Sinkhorn MPC with the number of Sinkhorn iterations $ S = 10,20,30 $ is $15.79, 12.04, 11.44 $, respectively.
As $ S $ increases, the accumulated cost approaches the cost for the unregularized case while the computational cost for obtaining $ P[k] $ grows.
Note that since we use MPC, the accumulated cost for the proposed method can be smaller than for the unregularized case as in this example.

\section{Conclusions}\label{sec:conclusion}
In this letter, we extended Sinkhorn MPC, which is a dynamical transport algorithm, to general input matrices.
Moreover, under the reachability of agents, we revealed its global convergence property, ultimate boundedness, and local asymptotic stability.
A numerical example validated the convergence of agents under Sinkhorn MPC close to the desired distribution and the cost-effectiveness of the proposed method.



\appendices

\section{Proof of Theorem~\ref{thm:convergence}}\label{app:lasalle}
	For notational simplicity, we drop the subscript $ P $.
        As a Lyapunov candidate function, we consider the entropic OT cost between $ x $ and $ x^\sfd $:
        \begin{equation}\label{eq:transport_cost_conti}
                \calE(x,x^\sfd) := \min_{P\in \calTone} \ \sum_{i,j\in \bbra{N}} C_{ij} (x) P_{ij} - \varepsilon \rmH(P) ,
        \end{equation}
        where $ C_{ij}(x) = c_{T_\rmh}^i (x_i,x_j^\sfd) = \|x_i - x_j^\sfd\|_{\scrG_i}^2,\ x^\sfd := [x_1^\sfd; \cdots ; x_N^\sfd] $.
        The time derivative of $ \calE(x(t),x^\sfd) $ along the trajectory of \eqref{eq:lq_mpc_dynamics} is given by
        \begin{align}
                \frac{\rmd}{\rmd t} \calE (x(t), x^\sfd) &= \nabla_x \calE (x(t), x^\sfd)^\top \dot{x} (t) \nonumber\\
    &= \sum_{i=1}^N \sum_{j=1}^N P_{ij}^* (x(t)) \nabla_{x_i} c_{T_\rmh}^i (x_i (t), x_j^\sfd)^\top \nonumber\\*
    &\quad \times \bar{A}_i \left(x_i(t) - x_i^\tmp (P^*(x(t))) \right)\nonumber\\
    &= 2\sum_{i=1}^N \sum_{j=1}^N P_{ij}^* (x(t))  (x_i(t) - x_j^\sfd)^\top \scrG_i \nonumber\\*
    &\quad \times \bar{A}_i \left(x_i(t) - x_i^\tmp (P^*(x(t))) \right) ,\nonumber
        \end{align}
        where $ \bar{A}_i := A_i - B_iB_i^\top \scrG_i $ and we used~\cite[Eq.~(9.6)]{Peyre2019}:
        \begin{equation}
                \nabla_{x_i} \calE(x,x^\sfd) = \sum_{j=1}^N P_{ij}^* (x) \nabla_{x_i} c_{T_\rmh}^i (x_i,x_j^\sfd) .
        \end{equation}
	Noting that $\sum_j P_{ij}^* (x(t)) x_j^\sfd = x_i^\tmp (P^*(x(t))) / N$ and $ P^*(x(t)) \in \calTone $, we obtain
  \begin{align*}
		&\frac{\rmd }{\rmd t} \calE (x(t), x^\sfd) = \frac{1}{N} \sum_i 2 \left(x_i(t) - x_i^\tmp (P^*(x(t)))\right)^\top \scrG_i  \\
		&\hspace{3.5cm} \times\bar{A}_i\left(x_i(t) - x_i^\tmp (P^*(x(t)))\right) \\
		&= \frac{1}{N} \sum_i \left(x_i(t) - x_i^\tmp (P^*(x(t)))\right)^\top \\
		&\qquad \times \scrG_i( \bar{A}_i\scrG_i^{-1} + \scrG_i^{-1}\bar{A}_i^\top )\scrG_i \left(x_i(t) - x_i^\tmp (P^*(x(t)))\right) .
  \end{align*}
	Here, $ \scrG_i^{-1} $ is known to satisfy the following Lyapunov equation\cite[Eq.~(8)]{Kleinman1970}:
		\[
			\bar{A}_i\scrG_i^{-1} + \scrG_i^{-1}\bar{A}_i^\top = - \Ee^{-A_iT_\rmh} B_iB_i^\top \Ee^{-A_i^\top T_\rmh} - B_iB_i^\top \preceq 0. 
		\]
  Then, it holds
\begin{align}\label{eq:lyap_decrease_energy}
	\frac{\rmd }{\rmd t} \calE (x(t), x^\sfd)
	\begin{cases}
		< 0, & x(t) \not\in \calM', \\
		= 0, & x(t) \in \calM'.
	\end{cases}
\end{align}

Lastly, by the unboundedness of the transport cost:
\begin{equation}\label{eq:unbounded}
        c_{T_\rmh}^i (x_i,x_j^\sfd) \rightarrow \infty \ {\rm as} \ \|x_i\| \rightarrow +\infty ,
\end{equation}
$ \calE $ is unbounded, and for any $ d\in \bbR $, the sublevel set $ \Omega_{\calE} (d) := \{x \in \bbR^{nN} : \calE(x,x^\sfd) \le d \}$ is bounded; see the proof of \cite[Theorem~1]{Ito2023sinkhorn}.
Then, by LaSalle's invariance principle~\cite[Theorem~4.4]{Khalil2002}, $ x(t) $ that follows \eqref{eq:lq_mpc_dynamics} starting in $ \Omega_{\calE}(d) $ converges to the largest invariant set in $ \calM' \cap \Omega_\calE (d) $. 
By the arbitrariness of $ d $ and the unboundedness of $ \calE $, we conclude that for any $ x(0) \in \bbR^{nN} $, $ x(t) $ converges to the largest invariant set in $ \calM' $ as $ t\rightarrow \infty $.
Lastly, $ x(t) \in \calM' $ yields \eqref{eq:invariant_system}.

\bibliographystyle{IEEEtran}

\bibliography{sinkhorn_cdc2023}

\end{document}